\documentclass{amsart}

\usepackage{amssymb}

\usepackage[all]{xy}

\newtheorem{thm}{Theorem}

\newtheorem{prop}[thm]{Proposition}

\theoremstyle{definition}
\newtheorem{defn}[thm]{Definition}

\newtheorem{say}[thm]{}

\newtheorem*{ack}{Acknowledgments}      

\newtheorem{defn-thm}[thm]{Definition--Theorem}  
\newtheorem{defn-lem}[thm]{Definition--Lemma}  

\theoremstyle{remark}

\renewcommand{\c}[0]{{\mathbb C}}  

\renewcommand{\o}[0]{{\mathcal O}} 
\newcommand{\z}[0]{{\mathbb Z}}
\newcommand{\n}[0]{{\mathbb N}}

\newcommand{\p}[0]{{\mathbb P}}

\newcommand{\q}[0]{{\mathbb Q}}

\newcommand{\qtq}[1]{\quad\mbox{#1}\quad}

\newcommand{\rank}[0]{\operatorname{rank}}
\newcommand{\mult}[0]{\operatorname{mult}}

\newcommand{\red}[0]{\operatorname{red}}    
\newcommand{\codim}[0]{\operatorname{codim}}    
    
\newcommand{\proj}[0]{\operatorname{Proj}}

\newcommand{\sing}[0]{\operatorname{Sing}}

\newcommand{\tsum}[0]{\textstyle{\sum}}

\def\into{\DOTSB\lhook\joinrel\to}

\def\loccoh#1.#2.#3.#4.{H^{#1}_{#2}(#3,#4)}

\DeclareMathAlphabet{\mathchanc}{OT1}{pzc}%
                                {m}{it}

\newcommand{\mdeg}[0]{\operatorname{Deg}}
\newcommand{\dres}[0]{\operatorname{\mathcal D\mathcal R}} 
\newcommand{\tprod}[0]{\textstyle{\prod}} 

\begin{document}
\bibliographystyle{amsalpha}

\title{Dual graphs of  exceptional divisors}
\author{J\'anos Koll\'ar}

\maketitle

Let $X$ be a complex algebraic or analytic variety.
Its local topology near  a point $x\in X$
is completely described by  the {\it link}  $L(x\in X)$,
which is obtained as  the intersection of
$X$ with a sphere of radius $0<\epsilon\ll 1$ centered at $x$.
A regular neighborhood of  $x\in X$ is 
homeomorphic to the cone over $L(x\in X)$; cf.\ \cite[p.41]{gm-book}.

One can study the local topology of
$X$ by  choosing a resolution of
 singularities  $\pi:Y\to X$ such that
$E_x:=\pi^{-1}(x)\subset Y$ is a simple normal crossing divisor
and then relating the topology of $E_x$ to the topology of the
 link  $L(x\in X)$.
This approach, initiated in \cite{mumf-top}, 
 has been especially successful for surfaces.

The topology of a simple normal crossing divisor
$E$ can in turn be understood in 2 steps.
First, the $E_i$ are smooth projective varieties,
and their topology is much studied.
A second layer of complexity comes from how the
components $E_i$ are glued together. This gluing process can be
naturally encoded by a finite cell  complex  ${\mathcal D}(E)$,
called the  {\it dual graph} or 
{\it dual  complex} of $E$;
see Definition \ref{snc-nc.sing.defn}. 
Given $x\in X$ and a  resolution  $\pi:Y\to X$, the 
 dual complex  ${\mathcal D}(E_x)$
depends on the resolution chosen, 
but its homotopy type does not; we 
denote it by $\dres(x\in X)$ 
(see \cite{MR2320738, MR2399025, 2011arXiv1102.4370A}).

Using this approach \cite{k-fg} proved that for every finitely presented
group $\Gamma$ there is a complex algebraic singularity
$(x\in X)$ (of dimension 3) 
such that  $\pi_1\bigl(L(x\in X)\bigr)\cong \Gamma$.
The proof starts with  first constructing 
a simple normal crossing variety $E$
such that
$\pi_1(E)\cong \pi_1\bigl({\mathcal D}(E)\bigr)\cong \Gamma$
and then realizing $E$ as the exceptional divisor on a resolution
of a singularity $(x\in X)$, yielding a chain of isomorphisms
$$
\pi_1\bigl(L(x\in X)\bigr)\cong\pi_1\bigl(\dres(0\in X)\bigr)\cong 
\pi_1\bigl({\mathcal D}(E)\bigr)\cong \Gamma.
$$

The aim of this note is to go further and prove that
not just $\pi_1\bigl(\dres(0\in X)\bigr)$ but $\dres(0\in X)$ can be arbitrary.

\begin{thm}\label{main.thm.1}
 Let $T$ be a  connected,  finite  cell complex.
Then there is a normal singularity $(0\in X_T)$
whose dual  complex $\dres(0\in X_T)$ is homotopy equivalent to $T$.
\end{thm}

It is interesting to connect properties of
$\dres(0\in X)$ with algebraic or geometric properties of the singularity
$(0\in X)$. 
A quasi projective variety $X$ has {\it rational singularities}
if for one (equivalently every) resolution of singularities $p:Y\to X$
and for every algebraic (or holomorphic) vector bundle $F$ on $X$,
the natural maps  $H^i(X, F)\to H^i(Y, p^*F)$ are isomorphisms.
That is, for purposes of computing cohomology of vector bundles,
$X$ behaves like a smooth variety.
See \cite[Sec.5.1]{km-book} for details.

It is known that if  $X$ has  rational singularities
then $\dres(0\in X)$ is  $\q$-acyclic, that is,
$H^i\bigl(\dres(0\in X), \q)=0$ for $i>0$, see for example
\cite[Lem.39]{k-fg}. 
The fundamental groups of the $\dres(0\in X) $
for rational singularities were determined in \cite{k-fg}:
these are exactly those finitely presented groups $G$
for which $H^1(G,\q)= H^2(G,\q)=0$
(sometimes called $\q$-superperfect groups).
Our next result determines the possible homotopy types of 
 $\dres(0\in X) $
for rational singularities.

\begin{thm}\label{main.thm.2}
 Let $T$ be a  connected, finite,  $\q$-acyclic cell complex.
Then there is a rational singularity $(0\in X)$
whose dual  complex $\dres(0\in X)$ is homotopy equivalent to $T$.
\end{thm}

While these are much stronger results than the fundamental group versions,
 most of the work needed to prove
these theorems 
was already done in \cite{k-t-lc-exmp, k-fg}.

The main technical result of \cite{k-fg} proves that for
every compact  simplicial complex $T$ there is a 
projective simple normal crossing variety $Z$ such that
${\mathcal D}(Z)$ is  homotopy equivalent to $T$,
while the main technical result of  \cite{k-t-lc-exmp} shows that
for every projective simple normal crossing variety $Z$
there is a normal singularity $(0\in X)$ and a 
{\em partial} resolution  $\pi:X'\to X$ such that
$Z\cong \pi^{-1}(0)\subset X'$.  

If $\dim Z\leq 4$ then $X'$ has only very simple singularities
which are easy to resolve. This was sufficient to control
$\pi_1\bigl(\dres(0\in X)\bigr)$. However, as $\dim Z$ increases,
$X'$ has more and more complicated singularities 
given locally by equations  of the form
$$
x_1\cdots x_n=t\cdot 
\det \left(
\begin{array}{ccc}
y_{11} & \cdots & y_{1m}\\
\vdots & & \vdots \\
y_{m1} & \cdots & y_{mm}
\end{array}
\right)
$$
where $n, m$ are arbitrary
and $(x_i, y_{ij}, t)$ are local coordinates.

If $\dim Z=2$, then the only singularity that appears is the
ordinary 3-fold double point
$(x_1x_2=ty_{11})$. The first somewhat complicated singularity  
$$
\bigl(x_1x_2x_3=t\cdot (y_{11}y_{22}-y_{12}y_{21})\bigr)\subset \c^8
$$
appears when $\dim Z=6$.
\medskip

In this paper we start with the varieties constructed in  
\cite{k-t-lc-exmp, k-fg} and
 resolve these singularities.
Surprisingly, the resolution process 
described in (\ref{indestup}--\ref{binres.say}) 
 leaves the  dual 
complex unchanged and we get the following.

\begin{thm}\label{main.thm.3} Let $Z$ be a projective  simple normal crossing
 variety of dimension $n$.
Then there is a normal singularity $(0\in X)$ of dimension $(n+1)$
and a resolution  $\pi:Y\to X$ such that
$E:=\pi^{-1}(0)\subset Y$ is a  simple normal crossing divisor 
and  its dual  complex  ${\mathcal D}(E)$
is naturally identified with   ${\mathcal D}(Z)$.
\end{thm}

\begin{say}[Open problems]\label{open.questions}{\ }

(\ref{open.questions}.1) It might be possible to describe 
all  complexes that occur on
resolutions of $(n+1)$-dimensional varieties. 
It is clear that $\dim {\mathcal D}(E)\leq n$
but I do not know any other restrictions.

Starting with an  $n$-dimensional  complex $T$, 
the constructions of \cite{k-t-lc-exmp, k-fg}
   give a $(2n+1)$-dimensional singularity $(0\in X)$
such that $\dres(0\in X)$ is homotopy equivalent to $T$. 
This increase of the dimension may not be necessary.

(\ref{open.questions}.2) The singularities constructed in Theorems
\ref{main.thm.1}--\ref{main.thm.2} are not isolated. It would be
interesting to construct isolated examples.

(\ref{open.questions}.3)   As shown by \cite{k-fg},  links of 
isolated  singularities are much more complicated topologically than
smooth projective varieties. As a starting point of further investigations,
 it would be useful
to understand the precise relationship between 
$\dres(0\in X)$ and the  topology  of the link of an 
isolated  singularity.

(\ref{open.questions}.4)   As we noted, given a 
singularity $(0\in X)$ and a resolution
$\pi:Y\to X$ such that $E:=\pi^{-1}(0)$ is a 
simple normal crossing divisor, the  homotopy type 
$\dres(0\in X)$ of 
the  dual  complex  ${\mathcal D}(E)$
does not depend on the choice of $\pi:Y\to X$. 

Note that if $p:X'\to X$ is a proper birational morphism
such that $E':=\pi^{-1}(0)$ is a 
simple normal crossing divisor, then 
the  dual  complex  ${\mathcal D}(E')$
is defined even if $X'$ is singular.

It is possible that
${\mathcal D}(E')$ is in fact homotopy equivalent to $\dres(0\in X)$
 as long as $X'$ has rational singularities. (The latter condition
is actually quite weak, for instance it  holds if none of the
strata of $E$  are contained in $\sing X$.)

(\ref{open.questions}.5)  Assume that $(X,\Delta)$ is dlt 
\cite[2.37]{km-book}. Since dlt implies rational,
$\dres(x\in X)$ a $\q$-acyclic for every point $x\in X$.
Furthermore, $\pi_1\bigl(\dres(x\in X)\bigr)=1$ by
\cite{k-shaf} and \cite{takayama}. Thus $\dres(x\in X)$ 
is contractible iff it is $\z$-acyclic. It would be very
interesting to decide  whether $\dres(x\in X)$ 
is contractible  or not.
For quotient singularities this is proved in 
\cite{2011arXiv1111.7177K}.

Related results are treated in \cite{k-ax} and \cite{hog-xu}.
\end{say}

\begin{defn}[Dual graphs]\label{snc-nc.sing.defn}
Let $X$ be a variety
with irreducible components $\{X_i: i\in I\}$.
We say that $X$ is a {\it simple normal crossing} variety
(abbreviated as {\it snc}) if the $X_i$ are smooth and 
 every point $p\in X$ has an open (Euclidean) neighborhood $p\in U_p\subset X$
and an embedding $U_p\into \c^{n+1}_p$ 
such that  $U_p\subset (z_1\cdots z_{n+1}=0)$.
A {\it stratum}
of $X$ is any irreducible component of an intersection
$\cap_{i\in J} X_i$ for some $J\subset I$.

The combinatorics of $X$ is encoded by a cell complex ${\mathcal D}(X)$
 whose vertices are labeled by the
 irreducible components of $X$ and for every stratum
$Z\subset \cap_{i\in J} X_i$ we attach a $(|J|-1)$-dimensional cell.
Note that for any $j\in J$ there is a unique  irreducible component
of $\cap_{i\in J\setminus\{j\}} X_i$ that contains $Z$;
this specifies the attaching map.
${\mathcal D}(X)$ is called the
{\it dual complex} or  {\it  dual graph} of $X$. 
(Although ${\mathcal D}(X)$ is not a simplicial complex in general, it is an
unordered $\Delta$-complex in the terminology of
\cite[p.534]{hatcher}.)
\end{defn}

\begin{defn}[Dual graphs associated to a singularity]\label{dres.defn}
 Let $X$ be a normal variety and
$x\in X$ a point. Choose a resolution of
 singularities  $\pi:Y\to X$ such that
$E_x:=\pi^{-1}(x)\subset Y$ is a simple normal crossing divisor.
Thus it has a dual  complex  ${\mathcal D}(E_x)$.

The dual graph of a normal surface singularity has a long history. Higher
dimensional versions appear in 
\cite{MR0506296, MR0466149, MR576865, friedman-etal}
but systematic investigations were started only recently; see
\cite{MR2320738, MR2399025, payne09, payne11}.

It is proved in \cite{MR2320738, MR2399025, 2011arXiv1102.4370A}
 that  the homotopy type of 
${\mathcal D}(E_x)$
is  independent of the resolution $Y\to X$. We denote it by
$\dres(x\in X)$. 
\end{defn}

 The proof of Theorem \ref{main.thm.1} starts with  the following,
which is a combination of Theorem 29 and Lemma 39
of \cite{k-fg}.

\begin{thm}\label{k-fg.thm}
 Let $T$ be a finite cell complex. Then there is a projective
simple normal crossing variety $Z_T$ such that
\begin{enumerate}
\item ${\mathcal D}(Z_T)$ is homotopy equivalent to $T$,
\item $\pi_1(Z_T)\cong \pi_1(T)$ and
\item $H^i\bigl(Z_T, \o_{Z_T}\bigr)\cong H^i(T, \c)$ for every $i\geq 0$.\qed
\end{enumerate}
\end{thm}

\begin{say}[Summary of the construction of \cite{k-t-lc-exmp}]
\label{summary.k-exmp}

Let $Z$ be a projective  local complete intersection variety
of dimension $n$ 
and choose any embedding  $Z\subset P$ into a smooth projective 
variety of dimension $N$.  
(We can take  $P=\p^N$ for $N\gg 1$.)
 Let $L$ be a sufficiently ample line bundle on $P$.
Let $Z\subset Y_1\subset P$ be the  complete intersection
of $(N-n-1)$  general sections of $L(-Z)$. Set
$$
Y:= B_{(-Z)}Y_1:=\proj_{Y_1}\tsum_{m=0}^{\infty}\o_{Y_1}(mZ).
$$
(Note that this is not the
 blow-up of $Z$ but  the
 blow-up of its inverse in the  class group.)

It is proved in \cite{k-t-lc-exmp} that 
the birational transform of $Z$ in $Y$ is  a Cartier divisor isomorphic to $Z$
and there is a  contraction morphism
$$
\begin{array}{ccc}
Z & \subset & Y\\
\downarrow && \hphantom{\pi}\downarrow\pi \\
0 & \in & X
\end{array}
\eqno{(\ref{summary.k-exmp}.1)}
$$
such that $Y\setminus Z\cong X\setminus \{0\}$.
If $Y$ is smooth then  $\dres(0\in X)={\mathcal D}(Z)$ and we are done
with  Theorem \ref{main.thm.1}.
However, the construction of \cite{k-t-lc-exmp} yields a smooth variety
$Y$  only if $\dim Z=1$ or $Z$ is smooth.
(It is easy to see that not every simple normal crossing variety $Z$ can be
realized as a hypersurface on a smooth variety, so this limitation is
not unexpected.)

Thus we need to understand the singularities of $Y$ and resolve them.

In order to do this, 
we need a very detailed description of the singularities of $Y$.
This is a local question, so we may assume that
 $Z\subset \c^N_{\mathbf x}$ is a complete intersection
  defined
by $f_1=\dots=f_{N-n}=0$. Let $Z\subset Y_1\subset \c^N$ be a general
complete intersection defined
by   equations
$$
h_{i,1} f_1+  \cdots  +h_{i,N-n} f_{N-n}  =  0
\qtq{for $i=1,\dots, N-n-1$. }
$$
Let $H=(h_{ij})$ be the $(N-n-1)\times (N-n)$ matrix of the system and 
$H_i$ the submatrix obtained by removing the $i$th column.
By \cite{k-t-lc-exmp} or \cite[Sec.3.2]{kk-singbook},
an open neighborhood of $Z\subset Y$ 
is defined by the equations
$$
\bigl(f_i=(-1)^i\cdot t\cdot \det H_i : i=1,\dots, N-n\bigr)
\subset \c^N_{\mathbf x}\times \c_t.
\eqno{(\ref{summary.k-exmp}.2)}
$$
Assume now that $Z$ has simple normal crossing singularities.
Up-to permuting the  $f_i$ and passing to a smaller open set, 
we may assume that
 $df_2,\dots, df_{N-n}$ are linearly independent everywhere along $Z$.
Then the singularities of $Y$ all come from the equation
$$
f_1=- t\cdot \det H_1.  
\eqno{(\ref{summary.k-exmp}.3)}
$$
Our aim is to write down  local normal forms for $Y$ along $Z$.

On $\c^N$ there is a stratification
$\c^N=R_0\supset R_1\supset \cdots$
where $R_i$ is the set of points where
$\rank H_1\leq (N-n-1)-i$. Since the $h_{ij}$ are general,
$\codim_W R_i=i^2$ and we may assume that every stratum of $Z$
is transversal to each $R_i\setminus R_{i+1}$ (\ref{detvar.say}).

Let $S\subset Z$ be any stratum and $p\in S$ a point
such that $p\in R_m\setminus R_{m+1}$.
We can choose local coordinates  $\{x_1,\dots, x_d\}$ and 
$\{y_{rs}: 1\leq r,s\leq m\}$  such that, in a neighborhood of $p$,
$$
f_1=x_1\cdots x_d\qtq{and} 
\det H_1=\det\bigl( y_{rs}: 1\leq r,s\leq m\bigr).
$$
Note that $m^2\leq \dim S=n-d$, thus we can add $n-d-m^2$ further
coordinates $y_{ij}$ to get a complete local coordinate system on $S$.

Then  the $n$ coordinates $\{x_k,y_{ij}\}$ determine a 
map
$$
\sigma:\c^N\times \c_t\to \c^n  \times \c_t
$$
such that $\sigma(Y)$ is defined by the equation
$$
x_1\cdots x_d=t\cdot \det\bigl( y_{rs}: 1\leq r,s\leq m\bigr).
$$
Since  $df_2,\dots, df_{N-n}$ are linearly independent along $Z$,
we see that   $\sigma|_Y$ is \'etale along $Z\subset Y$.
\end{say}

We can summarize these considerations as follows.

\begin{prop}\label{k-esmp.prop}
 Let $Z$ be a normal crossing variety of dimension $n$. 
Then there is
a normal singularity  $(0\in X)$ of dimension $n+1$ and a proper, birational
morphism $\pi:Y\to X$ such that $\red\pi^{-1}(0)\cong Z$
and  for every point $p\in \pi^{-1}(0)$
we can choose local  \'etale or analytic coordinates
 called 
$\{x_i: i\in I_p\}$ and $ \{y_{rs}:1\leq r, s\leq m_p\}$
(plus possibly other unnamed coordinates)
such that one can write the local equations of $Z\subset Y$ as 
$$
\bigl( \tprod_{i\in I_p}x_i=t=0\bigr)\subset \Bigl(\tprod_{i\in I_p}x_i=
t\cdot \det\bigl(y_{rs}:1\leq r, s\leq m_p\bigr)\Bigr)
\subset \c^{n+2}.\qed
$$
\end{prop}

\begin{say}[Determinantal varieties]\label{detvar.say}
We have used the following basic properties of determinantal varieties.
These are quite easy to prove directly; see \cite[12.2 and 14.16]{Harris95}
for a more general case.

Let $V$ be a smooth, affine variety,
and ${\mathcal L}\subset \o_V$
a finite dimensional sub vector space without common zeros.
Let  $H=\bigl(h_{ij}\bigr)$ be an $n\times n$ matrix
whose entries are general elements in ${\mathcal L}$.
For a point $p\in V$ set $m_p=\operatorname{corank} H(p)$. 
 Then there are local analytic coordinates
 $ \{y_{rs}:1\leq r, s\leq m_p\}$
(plus possibly other unnamed coordinates)
such that, in a neighborhood of $p$, 
$$
\det H= \det\bigl(y_{rs}:1\leq r, s\leq m_p\bigr).
$$
In particular, $\mult_p(\det H)=\operatorname{corank} H(p)$,
for every $m$  the set of points $R_m\subset V$ where
$\operatorname{corank} H(p)\geq m$ is a  subvariety
of pure codimension $m^2$ and $\sing R_m=R_{m+1}$.
\end{say}

\begin{say}[Inductive set-up for resolution]\label{indestup}
 The  object we try to resolve is  a triple
$$
(Y,E,F):=\bigl(Y, \tsum_{i\in I} E_i, \tsum_{j\in J} a_j F_j\bigr)
\eqno{(\ref{indestup}.1)}
$$
where $Y$ is a variety  over $\c$, $E_i, F_j$ are
codimension 1 subvarieties and $a_j\in \n$.
(The construction (\ref{summary.k-exmp}) produces a triple
$\bigl(Y, E:=Z, F:=\emptyset\bigr) $. The role of the $F_j$ is to keep track
of the exceptional divisors as we resolve the singularities of $Y$.)

We assume that   $E$ is a simple normal crossing variety and
for every point $p\in E$
there is a  (Euclidean) open neighborhood $p\in Y_p\subset Y$,
 an embedding $\sigma_p : Y_p\into \c^{\dim Y+1}$,
subsets $I_p\subset I$ and $J_p\subset J$, a natural number $m_p\in \n$
 and  local coordinates in $\c^{\dim Y+1}$ called 
$$
\{x_i: i\in I_p\},\ \{y_{rs}:1\leq r, s\leq m_p\},\  \{z_j: j\in J_p\}
\qtq{and} t
$$
(plus possibly other unnamed coordinates)
such that one can write the local equation of
 $\sigma_p(Y_p)\subset \c^{\dim Y+1}$ as
$$
 \tprod_{i\in I_p}x_i=
t\cdot \det\bigl(y_{rs}:1\leq r, s\leq m_p\bigr)\cdot \tprod_{j\in J_p}z_j^{a_j}.
\eqno{(\ref{indestup}.2)}
$$
Furthermore,   $\sigma_p(E_i)=(t=x_i=0)\cap \sigma_p(Y_p)$ for $i\in I_p$ and
$\sigma_p(F_j)=(z_j=0)\cap \sigma_p(Y_p)$ for $j\in J_p$.
(We do not impose any compatibility condition between the
local equations on overlapping charts.) 

We say that $(Y,E,F)$
is {\it resolved} at $p$ if $Y$ is smooth at $p$. 
\end{say}

The key technical result of the paper is the following.

\begin{prop}\label{main.res.prop} Let 
$(Y,E,F)$
be a triple as above.
Then there is a resolution of singularities
$\pi: \bigl(Y', E', F'\bigr)\to 
\bigl(Y, E, F\bigr)$
such that
\begin{enumerate}
\item $Y'$ is smooth and $E'$ is a 
simple normal crossing divisor,
\item $ E'=\pi^{-1}(E)$,
\item every stratum of $E'$ is mapped
birationally to a stratum of $E$ and
\item $\pi$ induces an identification
${\mathcal D}(E')={\mathcal D}(E)$.
\end{enumerate}
\end{prop}

Proof.  The resolution will be a composite of explicit blow-ups
of smooth subvarieties (except at the last step).
We use the local equations to describe the blow-up centers locally.
Thus we need to know which 
local subvarieties can be defined globally. 
For example,  choosing a
divisor $F_{j_1}$ specifies the local divisor $(z_{j_1}=0)$ at every point
$p\in F_{j_1}$.
Similarly, choosing two divisors $E_{i_1}, E_{i_2}$ gives
the local subvarieties   $(t=x_{i_1}=x_{i_2}=0)$  at every point
$p\in E_{i_1}\cap E_{i_2}$. (Here it is quite important that
the divisors $E_i$ are themselves smooth. The algorithm does not
seem to work if the $E_i$ have self-intersections.) Note that by contrast
$(x_{i_1}=x_{i_2}=0)\subset Y$ defines a local divisor 
which has no global meaning.
Similarly, the vanishing of any of the 
coordinate functions  $y_{rs}$ has no global meaning.

To a point $p\in \sing E$ we associate the local invariant
$$
\mdeg(p):=\bigl(\deg_x(p), \deg_y(p), \deg_z(p)\bigr)=
\bigl(|I_p|, m_p, \tsum_{j\in J_p} a_j\bigr).
$$
It is clear that $\deg_x(p)$ and $ \deg_z(p) $
do not depend on the local coordinates chosen. 
We see in (\ref{detres.say}) that $\deg_y(p)$ is also well defined
if $p\in \sing E$.
The degrees $\deg_x(p), \deg_y(p), \deg_z(p) $
are constructible and upper semi continuous functions on $\sing E$.

Note that $Y$ is smooth at $p$ iff either
$\mdeg (p)=(1,*,*)$ or $\mdeg (p)=(*,0,0)$. 
If $\deg_x(p)=1$ then we can rewrite the equation (\ref{indestup}.2)
as
$$
x'= t\cdot \tprod_{j}z_j^{a_j}
\qtq{where} 
x':=x_1+t\cdot \bigl(1-\det(y_{rs})\bigr)\cdot \tprod_{j}z_j^{a_j},
$$
so if $Y$is smooth then $\bigl(Y, E+ F\bigr)$
has only simple normal crossings along $E$. 
Thus the resolution constructed in Theorem \ref{main.thm.3}
is a log resolution.

The usual method of Hironaka would start by blowing up the {\em highest}
multiplicity points. This introduces new and rather complicated
exceptional divisors and I have not been able to understand
explicitly how the  dual complex changes.

In our case, it turns out to be  much better to
look at a  locus where $\deg_y(p)$ is maximal
but instead of maximizing $\deg_x(p)$ or $ \deg_z(p) $
we maximize the dimension. Thus we blow up 
subvarieties along which $Y$ is not equimultiple.
Usually this leads to a morass, but our equations separate
the variables into distinct groups which makes these blow-ups
easy to compute.

One can think of this as  mixing  the main step of the   Hironaka method
with the   order reduction for monomial ideals
(see, for instance, \cite[Step 3 of 3.111]{k-res}).

After some preliminary remarks about
blow-ups of simple normal crossing varieties 
 the proof of (\ref{main.res.prop})
is carried out in a series of steps (\ref{indestup}--\ref{binres.say}). 

We start with the locus where $\deg_y(p)$ is maximal and by a sequence
of blow-ups we eventually achieve that $\deg_y(p)\leq 1$
for every singular point  $p$. This, however, increases $\deg_z$.
Then in 3 similar steps we lower the maximum of $\deg_z$
until we achieve that $\deg_z(p)\leq 1$ for every singular point   $p$. 
Finally we take care of the  singular points where 
$\deg_y(p)+\deg_z(p)\geq 1$. \qed

\begin{say}[Blowing up simple normal crossing varieties]
\label{snc.blow-up.say}
Let $Z$ be a simple normal crossing variety and
$W\subset Z$ a  subvariety. 
We say that $W$ has {\it simple normal crossing} with $Z$ if
for each point $p\in Z$ there is an open neighborhood $Z_p$,
an embedding $Z_p\into \c^{n+1}$ and subsets
 $I_p, J_p\subset \{0, \dots, n\}$
such that
$$
Z_p=\bigl(\tprod_{i\in I_p} x_i=0\bigr)\qtq{and} 
W\cap Z_p=\bigl(x_j=0: j\in J_p\bigr).
$$
This implies that  
for every stratum $Z_{J}\subset Z$ the intersection
$W\cap Z_J$ is smooth (even scheme theoretically).  

If $W$ has  simple normal crossing with $Z$
then the blow-up $B_WZ$ is again a simple normal crossing variety.
If $W$ is one of the strata of $Z$, then 
${\mathcal D}(B_WZ)$ is obtained from ${\mathcal D}(Z)$
by removing the cell corresponding to $W$ and every other cell whose
closure contains it. Otherwise
${\mathcal D}(B_WZ)={\mathcal D}(Z)$.
(In the terminology of \cite[Sec.2.4]{kk-singbook},
  $B_WZ\to Z$ is a thrifty modification.)

As an example, let $Z=(x_1x_2x_3=0)\subset \c^3$. There are 7 strata
and  ${\mathcal D}(Z)$ is the 2-simplex whose vertices correspond to
the planes $(x_i=0)$.

Let us blow up a point $W=\{p\}\subset Z$ to
 get $B_pZ\subset B_p\c^3$. Note that the exceptional divisor 
$E\subset B_p\c^3$ is {\em not} a part of $B_pZ$
and $B_pZ$ still has 3 irreducible components.

If $p$ is the origin, then the triple intersection is removed
and  ${\mathcal D}(B_pZ)$ is the boundary of the  2-simplex.

If $p$ is not the origin, then $B_pZ$ still has 7 strata
naturally corresponding to the strata of $Z$ and
  ${\mathcal D}(B_pZ)$ is  the  2-simplex.

We will be interested in  situations where $Y$ is a hypersurface in $\c^{n+2}$
and  $Z\subset Y$ is a Cartier divisor that is a 
simple normal crossing variety.
Let $W\subset Y$ be a smooth, irreducible
 subvariety, not contained in $Z$  such that
\begin{enumerate}
\item the scheme theoretic intersection
$W\cap Z$ has simple normal crossing with $Z$
\item $\mult_{Z\cap W}Z=\mult_WY$. (Note that
this holds if $W\subset \sing Y$ and
$\mult_{Z\cap W}Z=2$.)
\end{enumerate}
Choose local coordinates  $(x_0,\dots, x_n,t)$ such that
$W=(x_0=\cdots x_i=0)$ and $Z=(t=0)\subset Y$.
Let $f(x_0,\dots, x_n,t)=0$ be the local equation of $Y$.

Blow up $W$ to get $\pi:B_WY\to Y$. 
Up to permuting the indices $0,\dots, i$, the blow-up $B_WY$ is covered
by coordinate charts described  by
the coordinate change
$$
\bigl(x_0, x_1, \dots, x_i, x_{i+1}, \dots, x_n,t\bigr)=
\big(x'_0, x'_1x'_0, \dots, x'_ix'_0, x_{i+1}, \dots, x_n,t\bigr).
$$
If $\mult_WY=d$ then the local equation of $B_WY$ in the above chart becomes
$$
(x'_0)^{-d}f\big(x'_0, x'_1x'_0, \dots, x'_ix'_0, x_{i+1}, \dots, x_n,t\bigr)=0.
$$
By assumption (2), 
$(x'_0)^{d} $ is also the largest power that divides
$$
f\big(x'_0, x'_1x'_0, \dots, x'_ix'_0, x_{i+1}, \dots, x_n,0\bigr),
$$
hence
$\pi^{-1}(Z)=B_{W\cap Z}Z$.

Observe finally that the conditions (1--2) can not be fulfilled in any
interesting way if $Y$ is smooth. Since we want $Z\cap W$ to be scheme
theoretically smooth,  if $Y$ is smooth then condition (1)
implies that $Z\cap W$ is disjoint from $\sing Z$.

(As an example, let $Y=\c^3$ and $Z=(xyz=0)$. Take $W:=(x=y=z)$. 
Note that $W$ is transversal to every irreducible component of $Z$
but $W\cap Z$ is a non-reduced point. The preimage of
$Z$ in $B_WY$ does not have  simple normal crossings.)

There are, however, plenty of  examples where $Y$ is singular
along $Z\cap W$ and these are exactly the singular points that
we want to resolve.

\end{say}

\begin{say}[Resolving the determinantal part]\label{detres.say}
Let $m$ be the largest size of a determinant occurring at a non-resolved point.
Assume that $m\geq 2$ and let
 $p\in Y$ be a non-resolved point with $m_p=m$. 

Away from $E \cup F $ 
the local equation of $Y$ is 
$$
 \tprod_{i\in I_p}x_i= \det\bigl(y_{rs}:1\leq r, s\leq m\bigr).
$$
Thus,
the singular set of 
$ Y_p\setminus (E \cup F)$ is
$$
\textstyle{\bigcup}_{(i, i')}
\bigl(\rank (y_{rs})\leq m-2\bigr)\cap 
 \bigl( x_{i}=x_{i'}=0\bigr)
$$
where the union runs through all 2-element subsets
$\{i,i'\}\subset I_p$. 
Thus the irreducible components of
$\sing Y\setminus (E \cup F)$
are in natural one-to-one correspondence with the
irreducible components of $\sing E$ and
the value of $m=\deg_y(p)$
is determined by the multiplicity of any of these
irreducible components at $p$.

Pick $i_1, i_2\in I$ and
we work locally with a subvariety 
$$
W'_p(i_1, i_2)
:=\bigl(\rank (y_{rs})\leq m-2\bigr)\cap\bigl(x_{i_1}=x_{i_2}=0\bigr).
$$
Note that $W'_p(i_1, i_2)$ is singular if $m>2$ and
the subset of its highest multiplicity points  is given by $\rank (y_{rs}) =0$.
Therefore the locally defined subvarieties
$$
W_p(i_1, i_2):=\bigl(y_{rs}=0 : 1\leq r,s\leq m\bigr)
\cap\bigl(x_{i_1}=x_{i_2}=0\bigr).
$$
glue together to  a well defined
global smooth subvariety $W:=W(i_1, i_2) $.

$E$ is defined by $(t=0)$ thus $E\cap W $ has the same
local equations as $W_p(i_1, i_2) $. In particular,
$E\cap W $ has simple normal crossings with $E$ and
$E\cap W $ is not a stratum of $E$;
its codimension in the stratum $(x_{i_1}=x_{i_2}=0) $ is $m^2$.

Furthermore, $E$ has multiplicity 2 along $E\cap W $,
hence (\ref{snc.blow-up.say}.2) also holds and so
$$
{\mathcal D}\bigl(B_{E\cap W}\bigr)= {\mathcal D}(E).
$$

We blow up $W \subset Y$. 
We will check that the new triple is again of the form
(\ref{indestup}).
The local degree $\mdeg(p)$
is unchanged over $Y\setminus  W $. The key assertion is that,
over $W $, the maximum value of $\mdeg(p)$
(with respect to the lexicographic ordering) decreases.
By repeating this procedure 
for every irreducible components of $\sing E$,
we  decrease the maximum value of
$\mdeg(p)$.
We can repeat this until we reach $\deg_y(p)\leq 1$ for every 
non-resolved point $p\in Y$.

(Note that this procedure requires an actual ordering of the
irreducible components of $\sing E$, which is a very non-canonical choice.
If a finite groups acts on $Y$, 
the resolution  usually can not be chosen  equivariant.)

Now to the local computation of the blow-up.
Fix a  point $p\in W $ and
set $I^*_p:=I_p\setminus\{i_1, i_2\}$.  
We write the local equation of $Y$  as
$$
x_{i_1}x_{i_2}\cdot L=t\cdot \det(y_{rs})\cdot R\qtq{where}
L:=\tprod_{i\in I^*_p}x_i\qtq{and}
R:= \tprod_{j\in J_p}z_j^{a_j}.
$$

There are two types of local charts on  the blow-up.

\begin{enumerate}
\item  There are two  charts of the first type. Up to interchanging
the subscripts $1,2$, these are given by the coordinate change
$$
(x_{i_1},x_{i_2},y_{rs} : 1\leq r,s\leq m)=
(x'_{i_1}, x'_{i_2}x'_{i_1}, y'_{rs}x'_{i_1} : 1\leq r,s\leq m).
$$ 
After setting $z_{w}:=x'_{i_1}$ the new local equation is
$$
 x'_{i_2}\cdot L=t\cdot  \det(y'_{rs})\cdot \bigl(z_{w}^{m^2-2} \cdot R\bigr).
$$  
The exceptional divisor is added to the $F$-divisors with coefficient
$m^2-2$
and the new degree is
$\bigl(\deg_x(p)-1, \deg_y(p), \deg_z(p)+m^2-2\bigr)$.

\item  There are $m^2$  charts of the second type. Up to 
re-indexing the  $m^2$ pairs $(r,s)$  these are given by the coordinate change
$$
(x_{i_1},x_{i_2},y_{rs} : 1\leq r,s\leq m)=
(x'_{i_1}y''_{mm},x'_{i_2}y''_{mm},y'_{rs}y''_{mm} : 1\leq r,s\leq m)
$$
except when $r=s=m$ where we set  $y_{mm}=y''_{mm} $.
It is convenient to set $y'_{mm} =1$ and
$z_w:=y''_{mm} $. Then the new local equation is
$$
x'_{i_1}x'_{i_2}\cdot L=t\cdot  \det\bigl(y'_{rs} : 1\leq r,s\leq m\bigr)
\cdot  \bigl(z_{w}^{m^2-2} \cdot R\bigr).
$$  
Note that  the $(m, m)$ entry of $(y'_{rs})$ is 1. 
By row and column operations we see that
$$
\det \bigl(y'_{rs}: 1\leq r,s,\leq m\bigr)=
\det \bigl(y'_{rs}-y'_{rm}y'_{ms}: 1\leq r,s,\leq m-1\bigr).
$$
By setting $y''_{rs}:=y'_{rs}-y'_{rm}y'_{ms}$ we have new local
equations
$$
x'_{i_1}x'_{i_2}L=t\cdot  
\det \bigl(y''_{rs}: 1\leq r,s,\leq m-1\bigr)
\cdot \bigl(z_{w}^{m^2-2}  \cdot R\bigr)
$$  
 and the new degree is
$\bigl(\deg_x(p), \deg_y(p)-1, \deg_z(p)+m^2-2\bigr)$.

\end{enumerate}
\medskip

{\it Outcome.} After these blow ups we have   a triple
$(Y,E,F)$
such that at non-resolved points the local equations are
$$
 \tprod_{i\in I_p}x_i=
t\cdot y\cdot \tprod_{j\in J_p}z_j^{a_j}
\qtq{or} 
 \tprod_{i\in I_p}x_i=
t\cdot \tprod_{j\in J_p}z_j^{a_j}.
\eqno{(\ref{detres.say}.3)}
$$
(Note that we can not just declare that $y$ is also a $z$-variable.
The $z_j$ are local equations of the divisors $F_j$ while
$(y=0)$ has no global meaning.)
\end{say}

\begin{say}[Resolving the monomial part]\label{monres.say}
Following  (\ref{detres.say}.3), the local equations are
$$
 \tprod_{i\in I_p}x_i=
t\cdot y^c\cdot \tprod_{j\in J_p}z_j^{a_j}\qtq{where $c\in \{0,1\}$.  }
$$
 We lower the degree 
of the $z$-monomial in 3 steps. 

{\it Step} 1. Assume that there is a non-resolved point with
$a_{j_1}\geq 2$. 

The singular set of $F_{j_1}$ is then
$$
\textstyle{\bigcup}_{(i, i')} \bigl(z_{j_1}=x_{i}=x_{i'}=0\bigr)
$$
where the union runs through all 2-element subsets
$\{i,i'\}\subset I$. 
 Pick an irreducible
component of it, call it 
$W(i_1, i_2, j_1):=\bigl(z_{j_1}=x_{i_1}=x_{i_2}=0\bigr) $.

Set $I^*_p:=I_p\setminus\{i_1, i_2\}$, $J^*_p:=J_p\setminus\{j_1\}$ and
 write the local equations as
$$
x_{i_1}x_{i_2}\cdot L=t z_j^{a_j}\cdot R\qtq{where}
L:=\tprod_{i\in I^*_p}x_i\qtq{and}
R:=y^c\cdot \tprod_{j\in J^*_p}z_j^{a_j}.
$$
There are 3 local charts on the blow-up:
\begin{enumerate}
\item $(x_{i_1},x_{i_2},z_j)=(x'_{i_1}, x'_{i_2}x'_{i_1}, z'_jx'_{i_1})$ and,
after setting $z_w:=x'_{i_1}$ 
the new local equation is
$$
 x'_{i_2}\cdot L=t\cdot  z_{w}^{a_j-2}{z'_j}^{a_j}\cdot R.
$$
The new degree is
$\bigl(\deg_x(p)-1, \deg_y(p), \deg_z(p)+a_j-2\bigr)$.
\item Same as above with the subscripts $1,2$ interchanged.
\item $(x_{i_1},x_{i_2},z_j)=(x'_{i_1}z'_j, x'_{i_2}z'_j, z'_j)$
with new local equation
$$
x'_{i_1}x'_{i_2}\cdot L=t \cdot {z'_j}^{a_j-2}\cdot R.
$$
The new degree is
$\bigl(\deg_x(p), \deg_y(p), \deg_z(p)-2\bigr)$.
\end{enumerate}

{\it Step} 2. Assume that there is a non-resolved point with
$a_{j_1}=a_{j_2}=1$. 

The singular set of $F_{j_1}\cap F_{j_2}$ is then
$$
\textstyle{\bigcup}_{(i, i')} \bigl(z_{j_1}=z_{j_2}=x_{i}=x_{i'}=0\bigr).
$$
where the union runs through all 2-element subsets
$\{i,i'\}\subset I$. 
 Pick an irreducible
component of it, call it 
$W(i_1, i_2, j_1, j_2):=\bigl(z_{j_1}=z_{j_2}=x_{i_1}=x_{i_2}=0\bigr) $.

Set $I^*_p:=I_p\setminus\{i_1, i_2\}$, $J^*_p:=J_p\setminus\{j_1, j_2\}$ and
we write the local equations as
$$
x_{i_1}x_{i_2}\cdot L=t z_{j_1}z_{j_2}\cdot R\qtq{where}
L:=\tprod_{i\in I^*_p}x_i\qtq{and}
R:=y^c\cdot \tprod_{j\in J^*_p}z_j^{a_j}.
$$
There are two types of local charts on the blow-up.

\begin{enumerate}
\item In the chart $(x_{i_1},x_{i_2},z_{j_1},z_{j_2})=
(x'_{i_1}, x'_{i_2}x'_{i_1}, z'_{j_1}x'_{i_1},z'_{j_2}x'_{i_1})$
 the new local equation is 
$$
 x'_{i_2}\cdot L=t\cdot  z'_{j_1}z'_{j_2}\cdot R.
$$ 
and the new degree is
$\bigl(\deg_x(p)-1, \deg_y(p), \deg_z(p)\bigr)$.
A similar chart is obtained by interchanging
the subscripts $i_1,i_2$.

\item In the chart $(x_{i_1},x_{i_2},z_{j_1},z_{j_2})=
(x'_{i_1}z'_{j_1}, x'_{i_2}z'_{j_1},z'_{j_1} ,z'_{j_2}z'_{j_1})$.
 the new local equation is 
$$
x'_{i_1}x'_{i_2}\cdot L=t \cdot z'_{j_2} \cdot R.
$$
The new degree is
$\bigl(\deg_x(p), \deg_y(p), \deg_z(p)-1\bigr)$.

A similar chart is obtained by interchanging
the subscripts $j_1,j_2$.
\end{enumerate}

By repeated application of these two steps we are reduced to the case
where $\deg_z(p)\leq 1$  at all  non-resolved points.

{\it Step} 3. Assume that there is a non-resolved point with
$\deg_y(p)=\deg_z(p)=1$. 

The singular set of $Y$ is 
$$
\textstyle{\bigcup}_{(i, i')} \bigl(y=z=x_{i}=x_{i'}=0\bigr).
$$
 Pick an irreducible
component of it, call it 
$W(i_1, i_2):=\bigl(y=z=x_{i_1}=x_{i_2}=0\bigr) $.
The blow up computation is the same as in Step 2.
\medskip

As before we see that at each step the conditions
(\ref{snc.blow-up.say}.1--2) hold, hence
${\mathcal D}(E)$ is unchanged. 
\medskip

{\it Outcome.} After these blow-ups we have   a triple
$(Y,E,F)$
such that at non-resolved points the local equations are
$$
 \tprod_{i\in I_p}x_i= t\cdot y, \quad 
 \tprod_{i\in I_p}x_i= t\cdot z_1 \qtq{or} 
 \tprod_{i\in I_p}x_i= t.
\eqno{(\ref{monres.say}.4)}
$$
As before, the $y$ and $z$ variables have different meaning,
but we can rename $z_1$ as $y$. Thus we have only one 
non-resolved local form left:
$\tprod x_i= t y $.
\end{say}

\begin{say}[Resolving the multiplicity 2 part]\label{binres.say}
 
Here we have a local equation   $x_{i_1}\cdots x_{i_d}=ty$
where $d\geq 2$. We would like to blow up
$(x_{i_1}=y=0)$, but, as we noted, 
this subvariety is not globally defined.
However,  a rare occurrence helps us out.
Usually the blow-up of a smooth subvariety
 determines its center uniquely. However, this is not the case
for codimension 1 centers. Thus we could get a globally well defined
blow-up even from centers that are not globally well defined.

Note that the inverse of $(x_{i_1}=y=0)$ in the local Picard 
group of $Y$ is $E_{i_1}=(x_{i_1}=t=0)$, which is globally defined.
Thus
$$
\proj_Y \tsum_{m\geq 0} \o_Y(mE_{i_1})
$$
is well defined, and locally it is isomorphic to the blow-up
 $B_{(x_{i_1}=y=0)}Y$. (A priori, we would need to take the
normalization of $B_{(x_{i_1}=y=0)}Y$, but it is actually normal.)
Thus we have 2 local charts.

\begin{enumerate}
\item $(x_{i_1},y)=(x'_{i_1}, y'x'_{i_1})$ and
the new local equation is
$\bigl(x_{i_2}\cdots x_{i_d}=ty'\bigr)$.
The new local degree is
$(d-1, 1,0)$.
\item $(x_{i_1},y)=(x'_{i_1}y', y')$ and
the new local equation is
$\bigl(x'_{i_1}\cdot x_{i_2}\cdots x_{i_d}=t\bigr)$.
The new local degree is
$(d, 0,0)$.
\end{enumerate}
\medskip

{\it Outcome.} After all these blow-ups we have   a triple
$\bigl(Y, \tsum_{i\in I} E_i, \tsum_{j\in J} a_j F_j\bigr)$
where  $\tsum_{i\in I} E_i $ is a simple normal crossing divisor and
$Y$ is smooth along   $\tsum_{i\in I} E_i $.
\medskip

This completes the proof of Proposition \ref{main.res.prop}. \qed
\end{say}

\begin{say}[Proof of Theorem \ref{main.thm.2}]
Assume that  $T$ is $\q$-acyclic.
Then, by (\ref{k-fg.thm}) there is a simple normal crossing variety
$Z_T$ such that $H^i\bigl(Z_T, \o_{Z_T}\bigr)=0$ for $i>0$.
Then \cite[Prop.9]{k-t-lc-exmp} shows that,
for $L$ sufficiently ample, 
 the 
singularity $(0\in X_T)$ constructed in (\ref{summary.k-exmp}) 
and (\ref{k-esmp.prop})
is rational. By (\ref{main.res.prop})
we conclude that $\dres(0\in X_T)\cong {\mathcal D}(Z_T)$
is homotopy equivalent to $T$.
\end{say}

 \begin{ack}
I thank M.~Kapovich, P.~Ozsv\'ath and S.~Payne for comments and corrections.
Partial financial support   was provided  by  the NSF under grant number 
DMS-07-58275.
\end{ack}


\def\cprime{$'$} \def\dbar{\leavevmode\hbox to 0pt{\hskip.2ex \accent"16\hss}d}
  \def\cprime{$'$} \def\cprime{$'$}
  \def\polhk#1{\setbox0=\hbox{#1}{\ooalign{\hidewidth
  \lower1.5ex\hbox{`}\hidewidth\crcr\unhbox0}}} \def\cprime{$'$}
  \def\cprime{$'$} \def\cprime{$'$} \def\cprime{$'$}
  \def\polhk#1{\setbox0=\hbox{#1}{\ooalign{\hidewidth
  \lower1.5ex\hbox{`}\hidewidth\crcr\unhbox0}}} \def\cdprime{$''$}
  \def\cprime{$'$} \def\cprime{$'$} \def\cprime{$'$} \def\cprime{$'$}
\providecommand{\bysame}{\leavevmode\hbox to3em{\hrulefill}\thinspace}
\providecommand{\MR}{\relax\ifhmode\unskip\space\fi MR }
\providecommand{\MRhref}[2]{%
  \href{http://www.ams.org/mathscinet-getitem?mr=#1}{#2}
}
\providecommand{\href}[2]{#2}

\noindent Princeton University, Princeton NJ 08544-1000

{\begin{verbatim}kollar@math.princeton.edu\end{verbatim}}

\end{document}